\documentclass[a4paper,12pt]{article}

\usepackage{amsmath,amsfonts,graphicx}

\usepackage{float}

\numberwithin{equation}{section}

\numberwithin{table}{section}

\title{Computational Experiments on $a^4+b^4+c^4+d^4=(a+b+c+d)^4$}

\author{Allan J. MacLeod\\Statistics, O.R. and Mathematics Group (Retired)\\
University of the West of Scotland,\\High St., Paisley,\\Scotland.   PA1 2BE\\
(e-mail: peediejenn@hotmail.com)}

\date{}

\begin{document}

\maketitle

\begin{abstract}
Computational approaches to finding non-trivial integer solutions of the equation in the title are discussed.
We summarize previous work and provide several new solutions.
\end{abstract}

\newpage

\section{Introduction}
The Diophantine equation
\begin{equation}\label{jmeq}
a^4+b^4+c^4+d^4=(a+b+c+d)^4
\end{equation}
was discussed by Jacobi and Madden in \cite{jm}, and has become known as the Jacobi-Madden equation.
They considered $a, b,c,d \in \mathbb{Z}$, but it is clear that the homogeneity of \eqref{jmeq}
means that we can consider rational values without loss of generality.

We have to assume at least three of the values are non-zero, because of Fermat's Last Theorem.
We, also, cannot have $a,b,c,d$ all of the same parity, so there must be a mixture of
positive and negative values.

The method used by Jacobi and Madden is based on the following simple, but remarkable, identity:
\begin{equation}\label{ident}
X^4+Y^4+(X+Y)^4=2(X^2+XY+Y^2)^2
\end{equation}

\eqref{jmeq} can be written
\begin{equation*}
a^4+b^4+(a+b)^4+c^4+d^4+(c+d)^4=(a+b)^4+(c+d)^4+(a+b+c+d)^4
\end{equation*}
and, using \eqref{ident}, we have
\begin{equation*}
(a^2+ab+b^2)^2+(c^2+cd+d^2)^2=((a+b)^2+((a+b)(c+d)+(c+d)^2)^2
\end{equation*}
if we ignore the common factor of $2$.

Let $F=a^2+ab+b^2$, $G=c^2+cd+d^2$ and $H=(a+b)^2+(a+b)(c+d)+(c+d)^2$, giving
\begin{equation*}
G^2=H^2-F^2=(H+F)(H-F)
\end{equation*}
so that
\begin{equation*}
\frac{H+F}{G}=\frac{G}{H-F}=t
\end{equation*}
where we will have $t \in \mathbb{Q}$. In fact, we have

\textbf{Lemma:} $t > 0$ for a non-trivial solution.

\textbf{Proof:} Each of $F, G, H$ is a variant of the basic quadratic form $Q(x,y)=x^2+xy+y^2$.

Defining $x=u+v, y=u-v$ gives $Q=3u^2+v^2 \ge 0$ which is only zero when $(u,v)=(0,0)=(x,y)$. Thus $F, G, H >0$, giving the result.

\vspace{0.5cm}

The first relation $H+F-Gt=0$ leads to the quadratic identity
\begin{equation*}
2a^2+3ab+2b^2+(a+b)(c+d)+(1-t)c^2+(2-t)cd+(1-t)d^2=0
\end{equation*}
which we can write in matrix-vector form as
\begin{equation}\label{mateq1}
(\begin{array}{llll}a&b&c&d\end{array}) \, \left( \begin{array}{rrrr}4&3&1&1\\3&4&1&1\\1&1&2(1-t)&2-t\\1&1&2-t&2(1-t)\end{array} \right) \, \left( \begin{array}{r}a\\b\\c\\d\end{array} \right) = 0
\end{equation}
where we have doubled the coefficients in the quadratic form to avoid fractions in the matrix. Call the $4 \times 4$ matrix $M_1$.

The relation $t(H-F)-G=0$ can be written, in a similar way, as
\begin{equation}\label{mateq2}
(\begin{array}{llll}a&b&c&d\end{array}) \, \left( \begin{array}{rrrr}0&t&t&t\\t&0&t&t\\t&t&2(t-1)&2t-1\\t&t&2t-1&2(t-1)\end{array} \right) \, \left( \begin{array}{r}a\\b\\c\\d\end{array} \right) = 0
\end{equation}
and we call this $4 \times 4$ matrix $M_2$.

Thus, the Jacobi-Madden problem reduces to finding a non-zero rational vector $\mathbf{v}$, with at least $3$ non-zero rational elements
and a non-zero rational $t$ such that
\begin{equation}\label{jmeqs}
\mathbf{v}^T \, M_1 \, \mathbf{v} = 0 = \mathbf{v}^T \, M_2 \, \mathbf{v}
\end{equation}

From $t=G/(H-F)$ we find
\begin{equation}\label{teq}
t=\frac{c^2+cd+d^2}{(a+c+d)(b+c+d)}
\end{equation}
which shows that a solution $(a,b,c,d)$ gives the same value of $t$ as $(a,b,d,c)$, $(b,a,c,d)$ and $(b,a,d,c)$. There are $24$ permutations
of a solution $(a,b,c,d)$, which are also solutions of the original problem, so they come in groups of $4$ giving $6$ different possible t-values.

For example, the solution found by Brudno $(5400,1770,-2634,955)$, which is used by Jacobi and Madden, leads to the t-values
$961/61$, $2521/325$, $1651/126$, $1777/1525$, $1423/1098$ and $511/450$. Note that for $t=961/61$ then $(t+1)/(t-1)=511/450$,
and the other $4$ t-values also form $\{t,(t+1)/(t-1)\}$ pairs. In fact, we have

\textbf{Lemma:} Given a non-trivial solution $(a,b,c,d)$ of \eqref{jmeq}, with $t$ given by \eqref{teq}, then
\begin{equation*}
\frac{t+1}{t-1}=\frac{a^2+ab+b^2}{(c+a+b)(d+a+b)}
\end{equation*}
the proof of which just involves a large amount of standard algebra, preferably done by a symbolic algebra package. Thus, we can assume
that, if $t=m/n$ with $m,n \in \mathbb{Z}$ and $\gcd(m,n)=1$, that $m$ and $n$ have opposite parities.

The present report discusses methods to compute other solutions, usually bigger in size. Since this problem could be of interest to amateurs,
I have tried to make the presentation as simple as possible.

\section{Quadric Intersection}
The first method is to use \eqref{mateq1} and \eqref{mateq2} directly. The intersection of two $4$-variable quadrics is fundamental to the
method of $4$-descent, used to find rational points on elliptic curves, see Merriman, Siksek and Smart \cite{mss} or the Ph.D. thesis of
Womack \cite{wom}.

\begin{table}[H]
\begin{center}
\caption{Solutions}
\begin{tabular}{lrrrr}
$\,$&$\,$&$\,$&$\,$&$\,$\\
t&a&b&c&d\\
$\,$&$\,$&$\,$&$\,$&$\,$\\
193/18&27385& 48150& 7590& -31764\\
511/450&-2634& 955& 5400& 1770\\
619/450&27385& -31764& 48150& 7590\\
1651/126&955& 5400& 1770& -2634\\
1141/666&7590& 27385& 48150& -31764\\
2041/150&-1229559& -1984340& 1022230& 107110\\
1423/1098&955& 1770& 5400& -2634
\end{tabular}
\end{center}
\end{table}

Mark Watkins of the Magma group in Sydney has an excellent preprint on the computational solution of such problems \cite{wat}. I have used
a Pari-GP version of this method for several years to compute points on hundreds of elliptic curves.

Applying this code using a fairly moderate search limit, with $t=m/n > 0$ and $m+n<3000$, gives the solutions in Table $2.1$.

There are only $3$ essentially different solutions in this Table. Increasing the search region but restricting to $m+n<499$
finds the extra solutions in Table $2.2$.

\begin{table}[H]
\begin{center}
\caption{Solutions}
\begin{tabular}{lrrrr}
$\,$&$\,$&$\,$&$\,$&$\,$\\
t&a&b&c&d\\
$\,$&$\,$&$\,$&$\,$&$\,$\\
31/6&53902630& 2542025& 35847220& -34122866\\
157/150&-841263& 792940& 44410& -3852350\\
181/150&-460945405& 189854902& 732896170& 303742360
\end{tabular}
\end{center}
\end{table}

The main problem with this method is that we do not know which values of $t$ to consider, so we start off by trying them all.
As we increase the search region, however, we need to restrict the choice of $t$.

We can reduce this by considering a change of variables used by Jacobi and Madden.

Let
\begin{equation}\label{tran1}
\left( \begin{array}{llll}a\\b\\c\\d\end{array} \right) = \left( \begin{array}{rrrr}0&0&2&2\\0&0&2&-2\\-1&-1&-1&0\\1&-1&-1&0\end{array} \right) \, \left( \begin{array}{r}p\\q\\r\\s\end{array} \right)
\end{equation}
and call the matrix in this transformation $C$, so that \eqref{jmeqs} become
\begin{equation}\label{jmeqs2}
\mathbf{w}^T \,C^T \, M_1 \, C \, \mathbf{w} = 0 = \mathbf{w}^T \, C^T \, M_2 \, C \, \mathbf{w}
\end{equation}
where $\mathbf{w}^T=(p,q,r,s)$.

Define $M_3=C^T \, M_1 \, C$ and $M_4=C^T \, M_2 \, C$ so that
\begin{equation}
M_3=\left( \begin{array}{rrrr}-2t&0&0&0\\-0&8-6t&-6t&0\\0&-6t&48-6t&0\\0&0&0&8t\end{array} \right) \hspace{0.5cm} M_4=\left( \begin{array}{rrrr}-2&0&0&0\\0&8t-6&-6&0\\0&-6&-6&0\\0&0&0&-8t\end{array} \right)
\end{equation}

Finally, define $M_5=(M_3-tM_4)/8$ giving
\begin{equation}
M_5=\left( \begin{array}{rrrr}0&0&0&0\\0&1-t^2&0&0\\0&0&6&0\\0&0&0&t^2+1\end{array} \right)
\end{equation}
so that we have that the quadric
\begin{equation}
(1-t^2)q^2+6r^2+(1+t^2)s^2=0
\end{equation}
must hold.

This clearly implies that $t^2>1$. We can also use the \textbf{Qfsolve} code from Denis Simon's ellrank package \cite{simon} to find
out whether this quadric has a rational solution for a specified value of $t$, rejecting those $t$ which have no solution.

\section{Quartic Equation}
In this section, we provide an alternative method of solution which also allows us to restrict greatly the
values of $t$ to be considered in lengthy computation. This was described by Tito Piezas III \cite{tp3} in a submission to the
\textbf{mathoverflow} web-site, where it elicited a very interesting response from Jeremy Rouse.

Let
\begin{equation}\label{tran2}
\left( \begin{array}{llll}a\\b\\c\\d\end{array} \right) = \left( \begin{array}{rrrr}1&-2&1&0\\1&-2&-1&0\\0&1&0&1\\0&1&0&-1\end{array} \right) \, \left( \begin{array}{r}p\\q\\r\\s\end{array} \right)
\end{equation}
with the matrix in this transformation called $D$, thus \eqref{jmeqs} becomes
\begin{equation}\label{jmeqs3}
\mathbf{w}^T \,D^T \, M_1 \, D \, \mathbf{w} = 0 = \mathbf{w}^T \, D^T \, M_2 \, D \, \mathbf{w}
\end{equation}

Define $M_{31}=D^T \, M_1 \, D$ and $M_{41}=D^T \, M_2 \, D$ so that
\begin{equation}
M_{31}=\left( \begin{array}{rrrr}14&-24&0&0\\-24&48-6t&0&0\\0&0&2&0\\0&0&0&-2t\end{array} \right) \hspace{1cm}
M_{41}=\left( \begin{array}{rrrr}2t&0&0&0\\0&-6&0&0\\0&0&-2t&0\\0&0&0&-2\end{array} \right)
\end{equation}

Next, define $M_{51}=(tM_{31}+M_{41})/2$ and $M_{61}=(tM_{41}-M_{31})/2$ giving
\begin{equation}
M_{51}=\left( \begin{array}{rrrr}8t&-12t&0&0\\-12t&-3(t^2-8t+1)&0&0\\0&0&0&0\\0&0&0&-(t^2+1)\end{array} \right)
\end{equation}
and
\begin{equation}
M_{61}=\left( \begin{array}{rrrr}t^2-7&12&0&0\\12&-24&0&0\\0&0&-(t^2+1)&0\\0&0&0&0\end{array} \right)
\end{equation}

In variable forms, we now have the quadrics
\begin{equation}
(t^2-7)p^2+24pq-24q^2=(t^2+1)r^2
\end{equation}
\begin{equation}
8tp^2-24tpq-3(t^2-8t+1)q^2=(t^2+1)s^2
\end{equation}

If we can find, for a given $t$, a solution $(p_0,q_0,r_0)$ $( q_0 \ne 0)$ to the first equation, we can parameterize using the
standard method. Simon's \textbf{Qfsolve} program tells us if the quadratic form is soluble and finds a solution, if possible.
This, then, becomes part of the sieving process for suitable $t$.

Let $x=p/q$ and $y=r/q$ so the first quadric is
\begin{equation*}
(t^2+1)y^2=(t^2-7)x^2+24x-24
\end{equation*}
with solution $x=x_0=p_0/q_0$ and $y=y_0=r_0/q_0$.

Then the line $y=y_0+k(x-x_0)$ will meet the quadric at one further point
\begin{equation*}
x=\frac{k^2x_0(t^2+1)-2ky_0(t^2+1)+t^2x_0-7x_0+24}{k^2(t^2+1)-t^2+7}
\end{equation*}
giving
\begin{equation}
\frac{p}{q}=\frac{k^2p_0(t^2+1)-2kr_0(t^2+1)+p_0(t^2-7)+24q_0}{q_0(k^2(t^2+1)-t^2+7)}
\end{equation}

Take the numerator for $p$ and the denominator for $q$, and substitute into
\begin{equation*}
(t^2+1)^2s^2=(t^2+1)(8tp^2-24tpq-3(t^2-8t+1)q^2)
\end{equation*}
and we have the quartic
\begin{equation}\label{quart}
Y^2=Ak^4+Bk^3+Ck^2+Dk+E
\end{equation}
where
\begin{equation*}
A=(t^2+1)^3(8p_0^2t-24p_0q_0t-3q_0^2(t^2-8t+1))
\end{equation*}
\begin{equation*}
B=16r_0t(t^2+1)^3(3q_0-2p_0)
\end{equation*}
\begin{equation*}
C=2(t^2+1)^2(8p_0^2t(t^2-7)+192p_0q_0t+3q_0^2(t^4-8t^3-6t^2-40t-7)+16r_0^2t(t^2+1))
\end{equation*}
\begin{equation*}
D=-16r_0t(t^2+1)^2(2p_0(t^2-7)+3q_0(t^2+9))
\end{equation*}
and
\begin{equation*}
E=(t^2+1)(8p_0^2t(t^2-7)^2+24p_0q_0t(t^2-7)(t^2+9)-
\end{equation*}
\begin{equation*}
3q_0^2(t^6-8t^5-13t^4-80t^3+35t^2-584t+49))
\end{equation*}

The quartic \eqref{quart} can be tested for local solubility - Simon includes Pari-GP code in \textbf{ellrank} - and those
which are not everywhere locally soluble can be rejected. It is also perfectly possible to reverse the order
in which the quadrics are considered. The smallest (in terms of $m+n$) t-values which give everywhere soluble quartics
are $31/6, 49/24, 67/42, \ldots$.

Using the quartic method (or something very similar), Seiji Tomita \cite{tom} found the following solutions.

\begin{table}[H]
\tiny
\begin{center}
\caption{Tomita solutions}
\begin{tabular}{lrrrr}
$\,$&$\,$&$\,$&$\,$&$\,$\\
t&a&b&c&d\\
$\,$&$\,$&$\,$&$\,$&$\,$\\
121/96&-889698809680& 687020381505&259448373800& 1526478290216 \\
121/96&22424373335225& 222795507072280&-237321095011880& 558974521862416 \\
181/150&802797814305& -626137906588& -150723250810& 1751113229630 \\
181/150&35966749745415& -360346958398438& 530920858665230& 377970149282480 \\
181/150&189854902& -460945405& 732896170& 303742360 \\
211/150&1229559&-1022230&1984340&-107110\\
211/150&561760&1493309&3597130&-1953890\\
373/150&-7929822455879583& 10830318289720550& 9309384955649330& 392431543415120\\
373/150&50627178820& 1357751663& 55867457830& -41572821650\\
709/450&1297734853& -1510410870& 500764020& 1768211850 \\
709/450&558360120& -701876813 & 753684930& 294589950 \\
3073/450&210240721& 396470430 & -336869940& 178944510\\
2851/1626&-2434795& 1945570& 1483582& 1858600\\
2977/2502&719130355& -2889516060& 4672341330& 2405612802
\end{tabular}
\end{center}
\end{table}

Studying the values of $t$ for the solutions, found so far, suggests the following

\textbf{Conjecture:} Let $t=m/n$ with $\gcd(m,n)=1$ and $m$ and $n$ of opposite parities.
If a solution exists, we will have $150|n$ or $25(6E+1)|(m-n)$, where $E \in \mathbb{Z}$. In the latter case
we have $6|n$.

I cannot believe I am the first person to think this! Can anyone prove or disprove this?

Using Simon's \textbf{Qfsolve} and \textbf{Qfparam} procedures, we can generate a multitude of quartics. I found that it was best
to apply Cremona's minimization and reduction methods \cite{cremred} to these quartics before searching for a point. With these methods, and a large
amount of computation the following new solutions were found.

\begin{table}[H]
\tiny
\begin{center}
\caption{New solutions}
\begin{tabular}{lrrrr}
$\,$&$\,$&$\,$&$\,$&$\,$\\
t&a&b&c&d\\
$\,$&$\,$&$\,$&$\,$&$\,$\\
499/474&3868630767650& 895775733285& 21271390911326& -4745425061560\\
511/150&-6714317914& 994485789915& -698106854980& 864417463190\\
3163/1350&-16515508578&10824551825&-15627586290&1711841340\\
18913/438&123140611690&446604426005&-96985017746&-25263498320
\end{tabular}
\end{center}
\end{table}

\section{Elliptic Curve}
Both the 4-descent and quartic methods have an underlying elliptic curve behind the problem.
To find this curve, we use the fact, from Merriman et al \cite{mss},
that a solution to \eqref{mateq1} and \eqref{mateq2} gives a point on the curve
\begin{equation*}
Y^2=\det(X\,M_1+M_2)
\end{equation*}
which can be given as
\begin{equation}\label{detq4}
Y^2=3t(7t-8)X^4-6(3t^3-7t+4)X^3-
\end{equation}
\begin{equation*}
3(t^4-8t^3+12t^2-7)X^2-6t(t^2-4t+3)X-3t^2
\end{equation*}

It is a standard fact, see chapter $3$ of Cremona  \cite{crem}, that a quartic with a rational point, is related to
the elliptic curve
\begin{equation}\label{q4toec}
y^2=x^3-27\,I \, x - 27\,J
\end{equation}
where $I$ and $J$ are the invariants of the quartic. The fundamental link is that rational $(X,Y)$ on \eqref{detq4} gets
mapped to a rational point with $x=3g_4(X)/4Y^2$ on \eqref{q4toec}, where $g_4$ is the quartic covariant of \eqref{detq4},
and $Y^2$ is given by \eqref{detq4}.

We find
\begin{equation}
I=9(t^8-16t^7+52t^6-48t^5+22t^4-176t^3+276t^2-144t+49)
\end{equation}
and
\begin{equation}
J=54K(t^8-16t^7+52t^6-144t^5+214t^4-176t^3+84t^2-48t+49)
\end{equation}
with $K=t^4-8t^3-6t^2+24t-7$.

Experiments with the right-hand-side of \eqref{q4toec} suggested it always factored, and it was reasonably straightforward
to find that $x=-9K$ gave $y=0$. Defining $z=x+9K$, and then $y=27v$ and $z=9u$ gives the fairly simple form
\begin{equation}\label{ecc}
E_t\,:\,v^2=u^3-3Ku^2+576t(t+1)(t-1)^3u
\end{equation}
Exactly the same elliptic curve comes from the quartic \eqref{quart} in the previous section. All the $p_0,q_0,r_0$ terms eventually vanish!

The elliptic curve $E_t$ has discriminant
\begin{equation}
\Delta=2^{16}\,3^6\,t^2(t+1)^2(t-1)^6(t^2+1)^2(t^4-16t^3+50t^2-80t+49)
\end{equation}
so $\Delta<0$ if $1.1742 < t < 12.483$ and $\Delta>0$ otherwise. If $\Delta<0$ the elliptic curve has one infinite component, whilst,
if $\Delta>0$, there is also a finite bounded component.
The curve is singular, for rational $t$, only for $|t|=1$ or $t=0$, but we saw in section $2$ that these values do not give solutions.

There is a clear rational point $u=0, \, v=0$ which is of order $2$. Numerical experiments suggest this is the only
finite torsion point, but there might well be specific values of $t$ giving extra torsion points.

These numerical experiments also suggested that the curve always has rank at least one. Results from ellrank indicated that
$u=48t$ gave a point, and it is easy to check that this gives $v=\pm 144t(t^2+1)$. If we double this point we find a point
where $u=4(t^2-2t-1)^2$.

Using ellrank and the Parity Conjecture, we find the ranks of the smallest t-values are given in Table $4.1$.

\begin{table}[H]
\begin{center}
\caption{Values of $t$}
\begin{tabular}{lr}
$\,$&$\,$\\
t&Estimated rank\\
$\,$&$\,$\\
31/6&2\\
49/24&2\\
67/42&1 or 3\\
79/54&2\\
97/72& 2 or 4\\
103/78&2 or 4\\
193/18&1 or 3
\end{tabular}
\end{center}
\end{table}
where we already have solutions for $t=31/6$ and $t=193/18$.

The basic fact about the rational points on an elliptic curve, over $\mathbb{Q}$, is that the points are finitely
generated. Thus, there exists a set of rational points $G_1,G_2,\ldots,G_r$ such that any rational point $P$ is such that
\begin{equation}\label{ecgen}
P=n_1G_1+n_2G_2+\ldots+n_rG_r+T
\end{equation}
where $n_1,\ldots,n_r \in \mathbb{Z}$ and $T$ is a torsion point.
$r$ is the rank of the elliptic curve and we assume $r \ge 1$ with $G_1=(48t,144t(t^2+1))$.

The elliptic curves \eqref{ecc} and \eqref{q4toec} can be easily transformed to one another. We have
\begin{equation*}
x=9(u-K)
\end{equation*}

Notice the direction of the relation of point on quartic to point on elliptic curve.
We \textbf{DO NOT} get a point on the quartic from every point on the elliptic curve.
In fact, I have never found a solution from $G_1$ or $2G_1$ with or without adding $(0,0)$. I wonder if there is a simple proof
of this? My attempts get bogged down in lots of variables.

Finding generators of elliptic curves is a highly non-trivial task. In fact, there is no known method guaranteed to work.
I initially used Magma's \textbf{TwoDescent} and \textbf{RationalPoints} procedures. Attempts to use Magma's
\textbf{FourDescent}, for large number of t-values, foundered as the computations take a long time,
admittedly on a not-very-fast machine. In March 2017, Pari
introduced the procedure \textbf{hyperellratpoints} which is an implementation of Michael Stoll's \textbf{ratpoints}. This meant that
I could use Pari for all the computations.

Given a set of generators, not necessarily of full rank, using
\eqref{ecgen} and the restriction $|n_i| \le L$, I generated points $P=(x,y)$ on \eqref{ecc}.
Then, I used Pari to factor
\begin{equation}
3g_4(X)-4\,x(P)\,Y^2
\end{equation}
to find a value of $X$ on \eqref{detq4} or other possible quartics.

For most acceptable t-values, we just find a single generator $G_1$.
For a few, we find a second generator, which may (or may not) lead to a solution of \eqref{jmeq}.
For $t=373/150$, we find $4$ generators with u-coordinates in the following Table.

\begin{table}[H]
\begin{center}
\caption{Generators for $t=373/150$}
\begin{tabular}{lr}
$\,$&$\,$\\
i&$u_i$\\
$\,$&$\,$\\
1&2984/25\\
2& 165858034880079528468553/154606810823279404439062500 \\
3&29529243840780598196578176/60686911309473227566225\\
4&184247616563459246903349991070216/16933216732179015462369769140625
\end{tabular}
\end{center}
\end{table}

Experiments show that the third generator must be included to give a solution of \eqref{jmeq}, so $|n_3| > 0$.
The numerical data all seem to suggest that solutions to \eqref{jmeq} all depend on one particular generator
being present in the expansion for a rational point.

This elliptic curve approach has found the following new reasonable-sized solutions

\begin{table}[H]
\tiny
\begin{center}
\caption{Elliptic Curve solutions}
\begin{tabular}{lrrrr}
$\,$&$\,$&$\,$&$\,$&$\,$\\
t&a&b&c&d\\
$\,$&$\,$&$\,$&$\,$&$\,$\\
1213/438&106185491830& 80795489585& 146163232960& -149806955726\\
1963/150&662971279500154&309770790508565&85290604949260& -371936154165950\\
1651/126&115711769730&58931380645&10424211666& -64829623500
\end{tabular}
\end{center}
\end{table}

For those values of $t$ given in Table $4.1$ without a solution, we looked at each value individually.
For $t=49/24$, the Birch and Swinnerton-Dyer conjecture
gives an estimate of the height of the other generator to be in the low hundreds, but within the computational
capabilities of my own software. By using the 2-isogenous curve, I found the second generator which gives the following
rather large solution
\begin{equation*}
a=-11590249845869269057824863556535439476779628603513075,
\end{equation*}
\begin{equation*}
b=12097338013880728917779953989473028810920897155225060,
\end{equation*}
\begin{equation*}
 c=3561881391291690403489592769705028154469958565069524,
 \end{equation*}
\begin{equation*}
 d=11315459134997579304238981942203181424806814023773640
\end{equation*}

For $t=79/54$, the 2-isogenous curve led nowhere, but the original curve finally gave up a second generator leading to
\begin{equation*}
a=246213540983698663206750 \hspace{1cm} b=4511618138222997480519985
\end{equation*}
\begin{equation*}
c= -4454458724579283498353610 \hspace{1cm} d=8579768155860334393439124
\end{equation*}

It is doubtful if these solutions could be found using either the quadric intersection or quartic-point methods.

\end{document}